\title{\vspace{-0.6cm} Nonnegative $k$-sums,
fractional covers, and probability of small deviations}
\author{
Noga Alon \thanks{Sackler School of Mathematics and Blavatnik
School of Computer Science, Tel Aviv University, Tel Aviv 69978,
Israel.
Email: {\tt nogaa@tau.ac.il}. Research supported in
part by an ERC Advanced grant and by a USA-Israeli BSF grant.}
\and
Hao Huang \thanks{Department of Mathematics, UCLA, Los
Angeles, CA, 90095. Email: {\tt huanghao@math.ucla.edu}.}
\and
Benny Sudakov\thanks{Department of Mathematics,
UCLA,  Los Angeles, CA 90095. Email: {\tt bsudakov@math.ucla.edu}.
Research supported in part by NSF CAREER award DMS-0812005 and
by a USA-Israeli BSF grant.}}

\documentclass[11pt]{article}
\usepackage{amsfonts,amssymb,amsmath,latexsym}
\usepackage{verbatim}

\date{}

\newtheorem{thm}{Theorem}[section]

\newtheorem{lemma}[thm]{Lemma}
\newtheorem{cor}[thm]{Corollary}

\newtheorem{conj}[thm]{Conjecture}

\newenvironment{pf}
      {\medskip\noindent{\bf Proof.}\hspace{1mm}}
      {\hfill$\Box$\medskip}

\def\qed{\ifvmode\mbox{ }\else\unskip\fi\hskip 1em plus 10fill$\Box$}

\begin{document}
\maketitle

\begin{abstract}
More than twenty years ago, Manickam, Mikl\'{o}s, and Singhi conjectured
that for any integers $n, k$ satisfying $n \geq 4k$, every set of $n$ real
numbers with nonnegative sum has at least $\binom{n-1}{k-1}$ $k$-element
subsets whose sum is also nonnegative.  In this paper we discuss the
connection of this problem with matchings and fractional covers of
hypergraphs, and with the question of estimating the probability
that the sum
of nonnegative independent random variables exceeds its expectation
by a given amount.
Using these connections together with some probabilistic techniques, we
verify the conjecture for $n \geq 33k^2$. This substantially improves the
best previously known exponential lower bound $n \geq e^{ck \log\log k}$.
In addition we prove a tight stability result showing that for
every $k$ and all sufficiently large $n$, every set of $n$ reals
with a nonnegative sum that does not contain a member whose sum
with any other $k-1$ members is nonnegative, contains at least 
$\binom{n-1}{k-1}+\binom{n-k-1}{k-1}-1$ subsets of cardinality
$k$ with nonnegative sum.
\end{abstract}

\section{Introduction}
\label{section_introduction}

Let $\{x_1,\cdots,x_n\}$ be a set of $n$ real numbers whose sum is
nonnegative. It is natural to ask the following question: how many subsets
of nonnegative sum must it always have?
The answer is quite straightforward, one can set $x_1=n-1$ and all the
other $x_i=-1$, which gives $2^{n-1}$ subsets. This construction is also
the smallest
possible since for every subset $A$, either $A$ or $[n]\backslash
A$ or both must have a nonnegative sum. Another natural question is, what
happens if we further restrict all the subsets to have a fixed size $k$?
The same example yields $\binom{n-1}{k-1}$ nonnegative $k$-sums consisting
of $n-1$ and $(k-1)$ $-1$'s. This construction is similar to the extremal
example in the Erd\H{o}s-Ko-Rado theorem \cite{erdos-ko-rado} which states
that for $n \geq 2k$, a family of subsets of size $k$ in $[n]$ with the
property that every two subsets have a nonempty intersection has size
at most $\binom{n-1}{k-1}$. However the relation between $k$-sum and
$k$-intersecting family is somewhat subtle and there is no obvious way
to translate one problem to the other.

Denote by $A(n,k)$ the minimum possible
number of nonnegative $k$-sums over
all possible choices of $n$ numbers $x_1, \cdots, x_n$
with $\sum_{i=1}^n x_i \geq 0$. For which values of $n$ and $k$,
is the construction $x_1=n-1, x_2=\cdots=x_n=-1$ best possible?
In other words, when can we guarantee that $A(n,k)=
\binom{n-1}{k-1}$?
This question was first raised by Bier and
Manickam \cite{bier, bier-manickam} in their study of
the so-called first distribution
invariant of the Johnson scheme. In 1987, Manickam
and Mikl\'{o}s \cite{manickam-miklos} proposed the
following conjecture, which in the language of the Johnson scheme
was also posed by Manickam and Singhi \cite{manickam-singhi} in 1988.

\begin{conj}
For all $n \geq 4k$, we have $A(n,k)=\binom{n-1}{k-1}$.
\end{conj}

In the Erd\H{o}s-Ko-Rado theorem, if $n<2k$, all the $k$-subsets form
an intersecting family of size $\binom{n}{k}>\binom{n-1}{k-1}$. But for
$n>2k$, the star structure, which always takes one fixed element
and $k-1$ other arbitrarily chosen elements, will do better than the set
of all $k$-subsets of the first $2k-1$ elements.
For a similar reason  we have the extra condition $n \geq 4k$
in the Manickam-Mik\'{o}s-Singhi conjecture. $\binom{n-1}{k-1}$ is not the
best construction when $n$ is very small compared to $k$. For example,
take
$n=3k+1$ numbers, $3$ of which are equal to $-(3k-2)$ and the other $3k-2$
numbers are $3$. It is easy to see that the sum is zero. On the other hand,
the nonnegative $k$-sums are those subsets consisting only of $3$'s,
which gives $\binom{3k-2}{k}$ nonnegative $k$-sums. It is not difficult to
verify that when $k>2$, $\binom{3k-2}{k} < \binom{(3k+1)-1}{k-1}$. However
this kind of construction does not exist for larger $n$.

The Manickam-Mik\'{o}s-Singhi conjecture has been open for more than two
decades. Only a few partial results of this conjecture are known so
far. The most important one among them is that the conjecture holds for
all $n$ divisible by $k$.  This claim can be proved directly by
considering a random partition of our set of numbers into
pairwise disjoint sets, each of size $k$, but it also
follows immediately from Baranyai's
partition theorem \cite{baranyai}. This theorem asserts
that if $k \mid n$, then the
family of all $k$-subsets of $[n]$ can be partitioned into disjoint
subfamilies so that each subfamily is a perfect $k$-matching. Since the
total sum is nonnegative, among the $n/k$ subsets  from each subfamily,
there must be at least one having a nonnegative sum. Hence there are no
less than $\binom{n}{k}/(n/k)=\binom{n-1}{k-1}$ nonnegative $k$-sums
in total. Besides this case, the conjecture is also known to be true
for small $k$. It is not hard to check it for $k=2$, and the case $k=3$
was settled by Manickam \cite{manickam}, and by Marino and Chiaselotti
\cite{marino-chiaselotti} independently.

Let $f(k)$ be the minimal number $N$ such that $A(n,k)=\binom{n-1}{k-1}$
for all $n \geq N$. The Manickam-Mikl\'{o}s-Singhi conjecture states that
$f(k) \leq 4k$. The existence of such function $f$ was first demonstrated
by Manickam and Mikl\'{o}s \cite{manickam-miklos} by showing $f(k)
\leq (k-1)(k^k+k^2)+k$.  Bhattacharya \cite{bhattacharya} found a new
and shorter proof of existence of $f$ later, but he didn't improve the
previous bound.  Very recently, Tyomkyn \cite{tyomkyn} obtained a better
upper bound $f(k)\leq k(4e \log k)^k \sim e^{ck\log\log k}$, which is
still exponential.

In this paper we discuss a connection between the
Manickam-Mikl\'{o}s-Singhi
conjecture and a problem about matchings in dense uniform hypergraph. We
call a hypergraph $H$ $r$-uniform if all the edges have size
$r$. Denote by $\nu(H)$ the matching number of $H$, which is the
maximum number of pairwise
disjoint edges in $H$. For our application, we need
the fact that if a $(k-1)$-uniform hypergraph on $n-1$ vertices has
matching number at most $n/k$, then its number of edges cannot exceed
$c\binom{n-1}{k-1}$ for some constant $c<1$ independent of $n,k$. This
is closely related to a special case of a long-standing open problem
of Erd\H{o}s \cite{erdos-matching}, who in 1965 asked to determine the
maximum  possible number of edges of an
$r$-uniform
hypergraph $H$ on $n$ vertices with matching number $\nu(H)$.
Erd\H{o}s conjectured that the
optimal case is when $H$ is a clique or the complement of a clique,
more precisely, for $\nu(H) < \lfloor n/r \rfloor$ the maximum
possible number of edges is given by the following  equation:
\begin{equation}
\max e(H) = \max \Bigg \{\binom{r[\nu(H)+1]-1}{r},~\binom{n}{r}
-\binom{n-\nu(H)}{r} \Bigg\}
\end{equation}

For our application to the Manickam-Mikl\'{o}s-Singhi conjecture, it
suffices to prove a weaker statement which bounds the number of edges
as a function of the fractional matching number $\nu^{*}(H)$ instead of
$\nu(H)$.  To attack the latter problem we combine duality
with a probabilistic technique
together with an inequality by Feige \cite{feige} which bounds the
probability that the sum of an arbitrary number of
nonnegative independent random variables exceeds its
expectation by a given amount.
Using this machinery, we obtain the first polynomial upper
bound $f(k) \leq 33k^2$, which substantially improves all the previous
exponential estimates.

\begin{thm}\label {mainthm_intro}
Given integers $n$ and $k$ satisfying $n\geq 33k^2$, for
any $n$ real numbers $\{x_1,\cdots,x_n\}$ whose sum is nonnegative,
there are at least $\binom{n-1}{k-1}$ nonnegative $k$-sums.
\end{thm}

Recall that earlier we mentioned the similarity between the 
Manickam-Mikl\'{o}s-Singhi conjecture and the Erd\H{o}s-Ko-Rado theorem.
When $n \geq 4k$, the conjectured extremal example is $x_1 = n-1, x_2=\cdots=x_n=-1$, where all the $\binom{n-1}{k-1}$ nonnegative $k$-sums
use $x_1$. For the Erd\H{o}s-Ko-Rado theorem when $n>2k$, the extremal family also consists of all the $\binom{n-1}{k-1}$ subsets containing one fixed element.
It is a natural question to ask if this kind of structure is forbidden, can we obtain a significant improvement on the $\binom{n-1}{k-1}$ bound?
A classical result of Hilton and Milner \cite{hilton-milner}
asserts that if $n>2k$ and no element is contained in every $k$-subset, 
then the intersecting family has size at most $\binom{n-1}{k-1}-
\binom{n-k-1}{k-1}+1$, with the extremal example being one of the 
following two.

\begin{itemize}
\item Fix $x \in [n]$ and $X \subset [n]\backslash \{x\}$, $|X|=k$. The family $\mathcal{F}_1$ consists of $X$ and all the $k$-subsets containing $x$ and intersecting with $X$.
\item Take $Y \subset [n]$, $|Y|=3$. The family $\mathcal{F}_2$ consists of all the $k$-subsets of
$[n]$ which intersects $Y$ with at least two elements.
\end{itemize}

It can be easily checked that both families are intersecting and $|\mathcal{F}_1| = \binom{n-1}{k-1}
-\binom{n-k-1}{k-1}+1$, $|\mathcal{F}_2| = 3 \binom{n-3}{k-2} + \binom{n-3}{k-3}$. When $k=3$,
$|\mathcal{F}_1|= |\mathcal{F}_2|$ and their structures are non-isomorphic.
For $k \ge 4$, $|\mathcal{F}_1| > |\mathcal{F}_2|$, so only the first construction is optimal.

Here we prove a Hilton-Milner type result 
about the minimum number of nonnegative $k$-sums.
Call a number $x_i$ \emph{large} if its sum with any other $k-1$
numbers $x_j$ is nonnegative. We prove that if no
$x_i$ is large, then the $\binom{n-1}{k-1}$ bound 
can be greatly improved. We also show that there are two extremal
structures, one of which is maximum for every $k$ and the other only for
$k=3$. This result can be considered as an analogue of the two extremal
cases mentioned above in the Hilton-Milner theorem.

\begin{thm} \label{hm}
For any fixed integer $k$ and sufficiently large $n$, and for any $n$ 
real numbers $x_1, \cdots, x_n$ with $\sum_{i=1}^n x_i \geq 0$, 
where no $x_i$ is large, the number $N$
of different nonnegative $k$-sums is at least 
$\binom{n-1}{k-1} + \binom{n-k-1}{k-1} -1$.
\end{thm}

For large $n$, Theorem \ref{hm} (whose statement is tight) 
improves the $\binom{n-1}{k-1}$ bound in the 
nonnegative $k$-sum problem to $\binom{n-1}{k-1}+\binom{n-k-1}{k-1}-1$ 
when large numbers are forbidden. This bound is asymptotically
$(2+o(1)) \binom{n-1}{k-1}$.

Call a number $x_i$  \emph{$(1-\delta)$-moderately large}, if there
are at least $(1-\delta)\binom{n-1}{k-1}$ nonnegative $k$-sums using
$x_i$, for some constant $0 \le \delta<1$. In particular, when $\delta=0$
this is the definition of a \emph{large} number. If there
is no $(1-\delta)$-moderately large number for some positive $\delta$,
we can prove a much stronger result asserting that at least a constant
proportion of the $\binom{n}{k}$ $k$-sums are nonnegative.
More precisely, we prove the
following statement.

\begin{thm} \label{moderate}
There exists a positive function $g(\delta, k)$, such that for any
fixed $k$ and $\delta$ and all sufficiently large $n$, the
following holds. For any set
of $n$ real numbers $x_1, \cdots, x_n$ with nonnegative sum
in which no member is
$(1-\delta)$-moderately large, the number $N$ of nonnegative $k$-sums
in the set
is at least $g(\delta, k) \binom{n}{k}$.
\end{thm}

The rest of this paper is organized as follows. In the next section
we present a quick proof of a slightly worse bound for the function
$f(k)$ defined above, namely,  we show that $f(k) \leq 2k^3$. The
proof
uses a simple estimate on the number of edges in a hypergraph with
a given matching number.  The proof of Theorem \ref{mainthm_intro}
appears in Section \ref{main}, where we improve this estimate using
more sophisticated probabilistic tools. In Section \ref{hi-mi} we prove the Hilton-Milner type results Theorem \ref{hm} and \ref{moderate}. The
last section contains some concluding remarks and open problems.

\section{Nonnegative $k$-sums and hypergraph matchings}\label{warmup}
In this section we  discuss the connection of the
Manickam-Mikl\'{o}s-Singhi conjecture and hypergraph matchings, and
verify the conjecture for $n \geq 2k^3$.

Without loss of generality, we can assume $\sum_{i=1}^n x_i=0$ and
$x_1 \geq x_2 \geq \cdots \geq x_n$ with $x_1>0$. If the sum of $x_1$
and the $k-1$ smallest numbers $x_{n-k+2},\cdots,x_n$ is nonnegative,
there are already $\binom{n-1}{k-1}$ nonnegative $k$-sums by taking $x_1$
and any other $k-1$ numbers. Consequently we can further assume that
$x_1+x_{n-k+2}+\cdots+x_n<0$. As all the numbers sum up to zero,
we have
\begin{equation}
x_2+\cdots+x_{n-k+1}>0
\end{equation}
Let $m$ be the largest integer not exceeding $n-k$ which
is divisible by $k$, then $n-2k+1 \leq m
\leq n-k$. Since the numbers are sorted in descending order, we have
\begin{equation}
x_2+\cdots+x_{m+1} \geq \dfrac{m}{n-k} (x_2+\cdots+x_{n-k+1})>0
\end{equation}
As mentioned in the introduction, the Manickam-Mik\'{o}s-Singhi conjecture
holds when $n$ is divisible by $k$ by Baranyai's partition theorem, thus
there are at least $\binom{m-1}{k-1} \geq \binom{n-2k}{k-1}$ nonnegative
$k$-sums using only numbers from $\{x_2,\cdots,x_{m+1}\}$. From now
on we are focusing on counting the number of nonnegative $k$-sums
involving $x_1$. If this number plus $\binom{n-2k}{k-1}$ is at least
$\binom{n-1}{k-1}$, then the Manickam-Mikl\'{o}s-Singhi conjecture is true.

Recall that in the proof of the
case $k \mid n$, if we regard all the negative
$k$-sums as edges in a $k$-uniform hypergraph, then the assumption that
all numbers add up to zero provides us the fact that this hypergraph
does not have a perfect $k$-matching. One can prove there are at least
$\binom{n-1}{k-1}$ edges in the complement of such a hypergraph, which
exactly tells the minimum number of nonnegative $k$-sums.
We utilize the same idea to estimate the number of nonnegative $k$-sums
involving $x_1$. Construct a $(k-1)$-uniform hypergraph $H$ on the
vertex set $\{2,\cdots,n\}$. The edge
set $E(H)$ consists of all the $(k-1)$-tuples $\{
i_1,\cdots,i_{k-1}\}$ corresponding to the negative $k$-sum
$x_1+x_{i_1}+\cdots+x_{i_{k-1}}<0$. Our goal is to show that $e(H)=|E(H)|$
cannot be too large, and therefore there must be lots of nonnegative
$k$-sums involving $x_1$.

Denote by $\nu(H)$ the matching number of our hypergraph $H$, which
is the maximum number of disjoint edges in $H$. By definition, every
edge corresponds to a $(k-1)$-sum which is less than $-x_1$, thus the
sum of the $(k-1)\nu(H)$ numbers corresponding to the vertices in the
maximal matching is less than $-\nu(H)x_1$. On the other hand, all
the remaining $n-1-(k-1)\nu(H)$ numbers are at most $x_1$. Therefore
$-x_1=x_2+\cdots+x_n \leq -\nu(H)x_1+(n-1-(k-1)\nu(H))x_1$. By solving
this inequality, we have the following lemma.

\begin{lemma}\label{matching}
The matching number $\nu(H)$ is at most $n/k$.
\end{lemma}

If the matching number of a hypergraph is known and $n$ is large with
respect to $k$, we are able to bound the number of its edges using
the following lemma. We denote by $\overline{H}$ the complement of
the hypergraph $H$.

\begin{lemma}\label{missingedges}
If $n>k^3$, any $(k-1)$-uniform hypergraph $H$ on $n-1$ vertices with
matching number at most $n/k$ has at least $\frac{1}{k+1}\binom{n-1}{k-1}$
edges missing from it.
\end{lemma}
\begin{pf}
Consider a random permutation $\sigma$ on the $n-1$ vertices
$v_1,\cdots,v_{n-1}$ of $H$. Let the random variables $Z_1=1$ if
$(\sigma(v_1),\cdots,\sigma(v_{k-1}))$ is an edge in $H$ and $0$
otherwise. Repeat the same process for the next $k-1$ indices and
so on.  Finally we will have at least $m\geq \frac{n-k}{k-1}$
random variables $Z_1,\cdots,Z_m$. Let $Z=Z_1+\cdots+Z_m$. $Z$
is always at most $n/k$ since there is no matching of size larger
than $n/k$. On the other hand, $\mathbb{E}Z_i$ is the probability
that $k-1$ randomly chosen vertices form an edge in $H$, therefore
$\mathbb{E}Z_i=e(H)/\binom{n-1}{k-1}$. Hence,
\begin{equation}\label{estimatematching}
\dfrac{n}{k} \geq \mathbb{E}Z = m \dfrac{e(H)}
{\binom{n-1}{k-1}}\geq \dfrac{n-k}{k-1} \dfrac{e(H)}{\binom{n-1}{k-1}}
\end{equation}
The number of edges missing is equal to
$e(\overline{H})=\binom{n-1}{k-1}-e(H)$. By \eqref{estimatematching},
$e(H) \leq (1-\frac{1}{k})\frac{n}{n-k}\binom{n-1}{k-1}$, therefore
\begin{eqnarray}
e(\overline{H}) &\geq&\Big[1-\Big(1-\dfrac{1}{k}\Big)
\dfrac{n}{n-k}\Big]\binom{n-1}{k-1} \nonumber\\
&\geq&\Big[1-\Big(1-\dfrac{1}{k}\Big)\dfrac{k^3}{k^3-k}\Big]
\binom{n-1}{k-1} \nonumber\\
&=&\dfrac{1}{k+1}\binom{n-1}{k-1}
\end{eqnarray}
\end{pf}

Now we can easily prove a polynomial upper bound for the function
$f(k)$ considered in the introduction, showing that
$f(k) \leq 2k^3$.
\begin{thm}\label{weakerresult}
If $n\geq2k^3$, then for any $n$ real numbers
$\{x_1,\cdots,x_n\}$ whose sum is nonnegative,
the number of nonnegative $k$-sums is at least $\binom{n-1}{k-1}$.
\end{thm}
\begin{pf}
By Lemma \ref{missingedges}, there are at least $\frac{1}{k+1}
\binom{n-1}{k-1}$ edges missing in $H$, which also gives a lower
bound for the number of nonnegative $k$-sums involving $x_1$.
Together with the previous $\binom{n-2k}{k-1}$ nonnegative $k$-sums
without using $x_1$,
there are at least $\frac{1}{k+1}\binom{n-1}{k-1}+\binom{n-2k}{k-1}$
nonnegative $k$-sums in total. We claim that this number is greater than
$\binom{n-1}{k-1}$ when $n\geq 2k^3$.  This statement is equivalent to
proving $\binom{n-2k}{k-1}/\binom{n-1}{k-1} \geq 1-1/(k+1)$, which can
be completed as follows:
\begin{eqnarray}
\binom{n-2k}{k-1}/\binom{n-1}{k-1} &=&\Big(1-\dfrac{2k-1}{n-1}\Big)
\Big(1-\dfrac{2k-1}{n-2}\Big)\cdots\Big(1-\dfrac{2k-1}{n-k+1}\Big)
\nonumber\\
&\geq& 1-\dfrac{(2k-1)(k-1)}{n-k+1} \nonumber\\
&\geq& 1-\dfrac{(2k-1)(k-1)}{2k^3-k+1} \nonumber\\
&\geq& 1-\dfrac{1}{k+1}
\end{eqnarray}
The last inequality is because $(2k-1)(k-1)(k+1)
=2k^3-k^2-2k+1 \leq 2k^3-k+1$.
\end{pf}

\section{Fractional covers and small deviations} \label{main}
The method above verifies the Manickam-Mikl\'{o}s-Singhi conjecture for
$n\geq 2k^3$ and improves the current best exponential lower bound $n\geq
k(4e \log k)^k$ by Tyomkyn \cite{tyomkyn}.  However if we look at Lemma
\ref{missingedges} attentively, there is still some room to improve
it. Recall our discussion of Erd\H{o}s' conjecture in the introduction:
if the conjecture is true in general, then in order to minimize the number
of edges in a $(k-1)$-hypergraph of a given matching number $\nu(H)=n/k$,
the hypergraph must be either a clique of size $(k-1)(n/k+1)-1$ or the
complement of a clique of size $n-n/k$.
\begin{equation}
e(H) \sim \min \bigg \{\binom{(1-1/k)n}{k-1},~\binom{n-1}{k-1}
-\binom{n-n/k}{k-1} \bigg\}
\sim (1-\dfrac{1}{k})^{k-1}\binom{n-1}{k-1} \leq \dfrac{1}{2}
\binom{n-1}{k-1}
\end{equation}
In this case, the number of edges missing from $H$ must be at least
$\frac{1}{2}\binom{n-1}{k-1}$, which is far larger than the bound
$\frac{1}{k+1}\binom{n-1}{k-1}$ in Lemma \ref{missingedges}. If in
our proof of Theorem \ref{weakerresult}, the coefficient before
$\binom{n-1}{k-1}$ can be changed to a constant instead of the
original $\frac{1}{k+1}$, the theorem can also be sharpened  to $n \geq
\Omega(k^2)$. Based on this idea, in the rest of this section we are going
to prove Lemma \ref{numberofedges}, which asserts that $e(H) \geq
c\binom{n-1}{k-1}$ for some constant $c$ independent of $n$ and $k$, and
can be regarded as a strengthening of Lemma \ref{missingedges}. Then we use
it to prove our main result of this paper, Theorem
\ref{maintheorem}. In order to
improve Lemma \ref{missingedges}, instead of using the usual matching
number $\nu(H)$, it suffices to consider its fractional relaxation,
which is defined as follows.
\begin{equation} \label{lp}
\begin{split}
\nu^{*}(H)&=\max \sum_{e \in E(H)}
w(e)~~~~~~~~~~~~w:E(H) \rightarrow [0,1]\\
&\textrm{subject to } \sum_{i \in e}w(e)
\leq 1~~~\textrm{for every vertex}~i.
\end{split}
\end{equation}

Note that $\nu^{*}(H)$ is always greater or equal than $\nu(H)$. On
the other hand for our hypergraph, we can prove the same upper
bound for the fractional matching number $\nu^{*}(H)$ as in Lemma
\ref{matching}. Recall that $H$ is the $(k-1)$-uniform hypergraph on
the $n-1$ vertices $\{2,\cdots,n\}$, whose edges are those $(k-1)$-tuples
$(i_1,\cdots,i_{k-1})$ corresponding to negative $k$-sums
$x_1+x_{i_1}+\cdots+x_{i_k}<0$.

\begin{lemma}\label{fractionalmatching}
The fractional matching number $\nu^{*}(H)$ is at most $n/k$.
\end{lemma}
\begin{pf}
Choose a weight function $w: E(H) \rightarrow [0,1]$ which optimizes
the linear program \eqref{lp} and gives the fractional matching
number $\nu^{*}(H)$, then $\nu^{*}(H)=\sum_{e \in E(H)} w(e)$. Two
observations can be easily made: (i) if $e \in E(H)$, then $\sum_{i \in
e} x_i<-x_1$; (ii) $x_i \leq x_1$ for any $i=2,\cdots,n$ since $\{x_i\}$
are in descending order. Therefore we can bound the fractional matching
number in a few steps.
\begin{eqnarray}
x_1+x_2+\cdots+x_n &= & x_1+\sum_{i=2}^{n} \Big(\sum_{i \in e} w(e)
\Big) x_i + \sum_{i=2}^{n} \Big(1-\sum_{i \in e} w(e)\Big) x_i \nonumber\\
&\leq& x_1+ \sum_{e \in E(H)} \Big(\sum_{i \in e} x_i\Big)w(e)
+ \sum_{i=2}^n \Big(1-\sum_{i \in e} w(e)\Big)x_1 \nonumber\\
&\leq& x_1 + \sum_{e \in E(H)} w(e) (-x_1) +
\Big(n-1-\sum_{e \in E(H)} \sum_{i \in e} w(e)\Big)x_1 \nonumber\\
&=&x_1-\nu^{*}(H)x_1+\Big(n-1-(k-1)\nu^{*}(H)\Big)x_1 \nonumber\\
&=&(n-k\nu^{*}(H))x_1
\end{eqnarray}
Lemma \ref{fractionalmatching} follows from this
inequality and our assumption that $x_1+\cdots+x_n=0$ and $x_1>0$.
\end{pf}

The determination of the fractional matching number is actually a linear
programming problem. Therefore we can consider its dual problem, which
gives the fractional covering number $\tau^{*}(H)$.

\begin{equation}\label{lp_covernumber}
\begin{split}
\tau^{*}(H)&=\min \sum_{i} v(i)~~~~~~~~~~~~v:V(H) \rightarrow [0,1]\\
&\textrm{subject to } \sum_{i \in e}v(i) \geq 1~~~\textrm{for every edge}~e.
\end{split}
\end{equation}

By duality we have $\tau^{*}(H)=\nu^{*}(H) \leq n/k$.  Getting a
upper bound for $e(H)$ is equivalent to finding a function $v:
V(H) \rightarrow [0,1]$ satisfying $\sum_{i \in V(H)} v(i) \leq n/k$
that maximizes the number of $(k-1)$-tuples $e$ where $\sum_{i \in e}
v(i) \geq 1$. Since this number is monotone increasing in every $v(i)$,
we can assume that it is maximized by a function $v$ with $\sum_{i \in
V(H)} v(i)=n/k$.

The following lemma was established by Feige \cite{feige}, and later
improved by He, Zhang, and Zhang \cite{he-zhang-zhang}. It bounds the
tail probability of the sum of independent nonnegative
random variables with given
expectations. It is stronger than Markov's inequality
in the sense that the
number of variables $m$ does not appear in the bound.

\begin{lemma}
\label{probabilistic_inequality}
Given $m$ independent nonnegative random variables
$X_1,\cdots,X_m$, each of expectation at most $1$, then
\begin{equation}
Pr\Big(\sum_{i=1}^m X_i <m+\delta \Big) \geq
\min \Big\{\delta/(1+\delta), \dfrac{1}{13}\Big\}
\end{equation}
\end{lemma}

Now we can show that the complement of the hypergraph $H$ has at least
constant edge density, which implies as a corollary that a constant
proportion of the $k$-sums involving $x_1$ must be nonnegative.

\begin{lemma} \label{numberofedges}
If $n \geq Ck^2$ with $C \geq 1$, and $H$ is a $(k-1)$-uniform hypergraph
on $n-1$ vertices with fractional covering number $\tau^{*}(H)=n/k$, then
there are at least $\Big(\dfrac{1}{13}-\dfrac{1}{2C}\Big)
\dfrac{(n-1)^{k-1}}{(k-1)!}$ $(k-1)$-sets which are not edges in $H$.
\end{lemma}

\begin{pf}
Choose a weight function $v: V(H) \rightarrow [0,1]$ which optimizes
the linear programming problem \eqref{lp_covernumber}.  Define a sequence
of $k-1$ independent and identically distributed random variables $X_1,
\cdots, X_{k-1}$, such that for any $1 \leq j \leq k-1, 2 \leq i \leq
n$, $X_j=v(i)$ with probability $1/(n-1)$. It is easy to compute the
expectation of $X_i$, which is
\begin{equation}
\mathbb{E}X_i=\dfrac{1}{n-1}\sum_{i=2}^{n} v(i)=\dfrac{n}{k(n-1)}
\end{equation}

Now we can estimate the number of $(k-1)$-tuples with sum less than
$1$. The probability of the event $\big\{\sum_{i=1}^{k-1}X_i<1\big\}$
is basically the same as the probability that a random $(k-1)$-tuple
has sum less than $1$, except that two random variables $X_i$ and $X_j$
might share a weight from the same vertex, which is not allowed for
forming an edge.  However, we assumed that $n$ is much larger than $k$,
so this error term is indeed negligible for our application. Note that
for $i_1<\cdots<i_{k-1}$, the probability that $X_j=v(i_j)$ for every
$1\leq j \leq k-1$ is equal to $1/(n-1)^{k-1}$.

\begin{eqnarray}
e(\overline{H})& =&\Big| \big\{i_1<\cdots<i_{k-1}:
v(i_1)+\cdots+v(i_{k-1})<1 \big \} \Big | \nonumber\\
&=& \dfrac{(n-1)^{k-1}}{(k-1)!} \sum_{\textrm{distinct~}i_1,
\cdots, i_{k-1}} Pr \Big[X_1=v(i_1),\cdots,X_{k-1}
=v(i_{k-1}); \sum_{i=1}^{k-1}X_i <1 \Big] \nonumber\\
& \geq& \dfrac{(n-1)^{k-1}}{(k-1)!}\bigg[Pr\Big(\sum_{i=1}^{k-1}
X_i<1\Big)- \sum_{l}\sum_{i\neq j}
Pr\Big(X_i=X_j=v(i_l)\Big)\bigg ] \nonumber\\
& \geq& \dfrac{(n-1)^{k-1}}{(k-1)!}\bigg [Pr\Big
(\sum_{i=1}^{k-1}X_i<1\Big)-
\dfrac{\binom{k-1}{2}}{n-1}\bigg ]  \label{sampling}\\
& \geq& \dfrac{(n-1)^{k-1}}{(k-1)!}\bigg [Pr\Big(\sum_{i=1}^{k-1}X_i
<1\Big)- \dfrac{1}{2C} \bigg ] \nonumber
\end{eqnarray}

The last inequality is because $n \geq Ck^2$ and $-1 \geq -3Ck+2C$ for
$C \geq 1, k \geq 1$, and the sum of these two inequalities implies that
$\dfrac{(k-1)(k-2)}{2(n-1)} \leq \dfrac{1}{2C}$.

Define $Y_i=X_i \cdot k(n-1)/n$ to normalize the expectations to
$\mathbb{E}Y_i=1$. $Y_i$'s are nonnegative because each vertex receives a
nonnegative weight in the linear program \eqref{lp_covernumber}. Applying
Lemma \ref{probabilistic_inequality} to $Y_1,\cdots,Y_{k-1}$ and
setting
$m=k-1$, $\delta=(n-k)/n$, we get
\begin{eqnarray}
Pr\Big(X_1+\cdots+X_{k-1} <1\Big)&=&Pr\Big(Y_1
+\cdots+Y_{k-1} <k(n-1)/n\Big) \nonumber \\
&\geq &\min \Big\{\dfrac{n-k}{2n-k},
\dfrac{1}{13}\Big\} \label{inequality_forY}
\end{eqnarray}
When $n>Ck^2$ and $k \geq 2$, $C \geq 1$, we have
\begin{equation}
\dfrac{n-k}{2n-k}>\dfrac{Ck^2-k}
{2Ck^2-k}=\dfrac{Ck-1}{2Ck-1} \geq \dfrac{1}{13}
\end{equation}

Combining \eqref{sampling} and \eqref{inequality_forY}
we immediately obtain Lemma \ref{numberofedges}.
\end{pf}

\begin{cor} \label{cor_numberofedges}
If $n \geq Ck^2$ with $C \geq 1$, then there are at least
$\Big(\dfrac{1}{13}-\dfrac{1}{2C}\Big)\dfrac{(n-1)^{k-1}}{(k-1)!}$
nonnegative $k$-sums involving $x_1$.
\end{cor}

Now we are ready to prove our main theorem:

\begin{thm}\label{maintheorem}
If $n \geq 33k^2$, then for any $n$ real numbers $x_1,\cdots,x_n$
with $\sum_{i=1}^n x_i \geq 0$,
the number of nonnegative $k$-sums is at least $\binom{n-1}{k-1}$.
\end{thm}
\begin{pf}
By the previous discussion, we know that there are at least
$\binom{n-2k}{k-1}$ nonnegative $k$-sums using only $x_2,
\cdots,x_n$. By Corollary \ref{cor_numberofedges},
there are at least $\Big(\dfrac{1}{13}-\dfrac{1}{2 \cdot
33}\Big)\dfrac{(n-1)^{k-1}}{(k-1)!} \geq \dfrac{2}{33}
\dfrac{(n-1)^{k-1}}{(k-1)!}$ nonnegative $k$-sums involving $x_1$. In
order to prove the theorem, we only need to show that for $n \geq 33k^2$,
\begin{equation}
\dfrac{2}{33} \dfrac{(n-1)^{k-1}}{(k-1)!}+
\binom{n-2k}{k-1} \geq \binom{n-1}{k-1}
\end{equation}
Define an infinitely differentiable function $g(x)
=\dbinom{x}{k-1}=\dfrac{x(x-1)\cdots(x-k+2)}{(k-1)!}$.
It is not difficult to see $g''(x)>0$ when $x>k-1$. Therefore
\begin{equation}\label{newtonmethod}
\binom{n-1}{k-1}-\binom{n-2k}{k-1}=g(n-1)-g(n-2k)
\leq [(n-1)-(n-2k)]g'(n-1)=(2k-1)g'(n-1)\\
\end{equation}
\begin{eqnarray}
g'(x)&=&\dfrac{(x-1)(x-2)\cdots(x-k+2)}{(k-1)!}+\dfrac{x(x-2)
\cdots(x-k+2)}{(k-1)!}+\cdots+\dfrac{x(x-1)\cdots(x-k+3)}{(k-1)!}
\nonumber\\
&\leq& (k-1) \dfrac{x(x-1)\cdots(x-k+3)}{(k-1)!} \nonumber\\
&\leq&  (k-1) \dfrac{x^{k-2}}{(k-1)!} \label{derivative}
\end{eqnarray}
Combining \eqref{newtonmethod} and \eqref{derivative},
\begin{equation}
\binom{n-1}{k-1}-\binom{n-2k}{k-1} \leq (2k-1)f'(n-1) \leq
(2k-1)(k-1) \dfrac{(n-1)^{k-2}}{(k-1)!}
\leq \dfrac{2}{33}\dfrac{(n-1)^{k-1}}{(k-1)!}
\end{equation}
The last inequality follows from our assumption $n \geq 33k^2$.
\end{pf}

\section{Hilton-Milner type results} \label{hi-mi}
In this section we prove two Hilton-Milner type results about the 
minimum number of nonnegative $k$-sums. The first theorem asserts that 
if $\sum_{i=1}^n x_i \geq 0$ and no $x_i$ is large, then there are at least  $\binom{n-1}{k-1} + \binom{n-k-1}{k-1}-1$ nonnegative 
$k$-sums.\\
\\
\noindent \textbf{Proof of Theorem \ref{hm}.}
We again assume that $x_1 \geq \cdots \geq x_n$ and $\sum_{i=1}^n x_i$ is zero. Since $x_1$ is not large, we know that
there exists a $(k-1)$-subset $S_1$ not containing $1$, such that $x_1 + \sum_{i \in S_1} x_i <0$. Suppose $t$ is the largest
integer so that there are $t$ subsets $S_1, \cdots, S_t$, such that for any $1 \leq j \leq t$, $S_j$ is 
disjoint from $\{1, \cdots, j\}$, has size at most $j(k-1)$ and
$$x_1 + \cdots +x_j + \sum_{i \in S_j} x_i < 0.$$

As we explained above $t \geq 1$ and since $x_1 \geq \cdots \geq x_n$ we may also assume that $S_j$ consists of the last $|S_j|$ indices in $\{1, \cdots, n\}$.
By Corollary \ref{cor_numberofedges}, for sufficiently large $n$, 
there are at least $\frac{1}{14} \binom{n-1}{k-1}$ nonnegative $k$-sums using $x_1$.
Note also that after deleting $x_1$ and $\{x_i\}_{i \in S_1}$, the sum of the remaining $n-1-|S_1| \geq n-k$ numbers is nonnegative. Therefore, again by Corollary \ref{cor_numberofedges},
there are at least $\frac{1}{14} \binom{n-k-1}{k-1}$ nonnegative $k$-sums using $x_2$ but not $x_1$. In the next step,
we delete $x_1, x_2$ and $\{x_i\}_{i \in S_2}$ and bound the number of nonnegative $k$-sums involving $x_3$ but neither $x_1$ or $x_2$
by $\frac{1}{14} \binom{n-2k-1}{k-1}$. Repeating this process for $30$ steps, we obtain
\begin{eqnarray}
N &\geq & \dfrac{1}{14} \left[ \binom{n-1}{k-1} + \binom{n-k-1}{k-1} + \cdots + \binom{n-29k-1}{k-1}\right] \nonumber \\
&> &  \dfrac{30}{14} \binom{n-29k-1}{k-1} > 2 \binom{n-1}{k-1} > \binom{n-1}{k-1} + \binom{n-k-1}{k-1} -1  \nonumber
\end{eqnarray}
where here we used the fact that
 $\frac{30}{14}>2$ and $n$ is sufficiently large (as a function of
$k$).

If $2 \leq t < 30$, by the maximality of $t$, we know that the sum of $x_{t+1}$ with any $(k-1)$ numbers with indices not in $\{1, \cdots, t+1\} \cup S_t$ is nonnegative. This
gives us $\binom{n-(t+1)-|S_t|}{k-1} \geq \binom{n-tk-1}{k-1}$ nonnegative $k$-sums. We can also replace $x_{t+1}$ by any $x_i$ where $1 \leq i \leq t$ and the new $k$-sum is still nonnegative since $x_i \geq x_{t+1}$.
Therefore,
$$N \geq (t+1) \binom{n-tk-1}{k-1} \geq (t+1) \binom{n-29k-1}{k-1} > 
2 \binom{n-1}{k-1}$$ for sufficiently large $n$. 
Thus the only remaining case is $t=1$.

Recall that $x_1$ is not large, and hence $x_1 + (x_{n-k+2} + \cdots + x_n)
< 0$.  Suppose $I$ is a $(k-1)$-subset of $[2,n]$ such that $x_1 +
\sum_{i \in I} x_i < 0$. If $2 \in I$, then $x_1 + x_2 + \sum_{i \in
I \backslash \{2\}} x_i < 0$, this contradicts the assumption $t=1$
since $|I \backslash\{2\}|=k-2 \leq 2(k-1)$.  Hence we can assume that
all the $(k-1)$-subsets  $I_1, \cdots, I_m$ corresponding to negative
$k$-sums involving $x_1$ belong to the interval $[3,n]$. Let $N_1$ be
the number of nonnegative $k$-sums involving $x_1$, and let $N_2$ be the
number of nonnegative $k$-sums using $x_2$ but not $x_1$, then
$$
N \geq N_1 + N_2 = \left[\binom{n-1}{k-1}-m \right] + N_2.
$$
In order to prove $N \geq \binom{n-1}{k-1}+\binom{n-k-1}{k-1}-1$, we only need to establish the following inequality
\begin{equation}\label{counting}
N_2 \geq \binom{n-k-1}{k-1} + m - 1.
\end{equation}
Observe that the subsets $I_1, \cdots, I_m$ satisfy some additional properties. First of all, if two sets $I_i$ and $I_j$ are disjoint, then by definition, $x_1 + \sum_{l \in I_i} x_l <0$ and
$x_2 + \sum_{l \in I_j} x_l \leq  x_1 + \sum_{l \in I_j} x_l<0$, summing them up gives
$x_1 + x_2 + \sum_{l \in I_i \cup I_j} x_l < 0$ with $|I_i \cup I_j|=2(k-1)$,
which again contradicts the assumption $t=1$. Therefore we might
assume that $\{I_i\}_{1 \leq i \leq m}$ is an intersecting family.
By the Erd\H{o}s-Ko-Rado theorem,
$$m \leq \binom{(n-2)-1}{(k-1)-1}=\binom{n-3}{k-2}.$$
The second observation is that if a $(k-1)$-subset $I \subset [3,n]$ 
is disjoint from some $I_i$, then $x_2 + \sum_{i \in I} x_i \geq 0$. 
Otherwise if $x_2 + \sum_{i \in I} x_i < 0$ and
$x_1 + \sum_{k \in I_i} x_k < 0$, for the same reason this contradicts
$t=1$. Hence $N_2$ is bounded from below by the number of $(k-1)$-subsets
$I \subset [3,n]$ such that $I$ is disjoint from at least one of $I_1,
\cdots, I_m$. Equivalently we need to count the distinct $(k-1)$-subsets
contained in some $J_i = [3,n] \backslash I_i$, all of which have sizes
$n-k-1$. By the real version of the Kruskal-Katona theorem (see Ex.13.31(b)
in \cite{lovasz}), if $m= \binom{x}{n-k-1}$ for some positive real number
$x \geq n-k-1$, then $N_2 \geq \binom{x}{k-1}.$ On the other hand, it is
already known that $1 \leq m \leq \binom{n-3}{k-2} = \binom{n-3}{n-k-1}$,
thus $n-k-1 \leq x \leq n-3$. The only remaining step is to verify the
following inequality for $x$ in this range,
\begin{equation} \label{kk}
\binom{x}{k-1} \geq \binom{n-k-1}{k-1} + \binom{x}{n-k-1} -1.
\end{equation}
Let $f(x) = \binom{x}{k-1} - \binom{x}{n-k-1}$, note that when $x \leq n-4 = (k-2)+(n-k-2)$,
\begin{eqnarray}
f(x+1)-f(x) &=& \left[\binom{x+1}{k-1} - \binom{x+1}{n-k-1}\right] - \left[\binom{x}{k-1} - \binom{x}{n-k-1}\right] \nonumber\\
&=& \binom{x}{k-2} - \binom{x}{n-k-2} \geq 0 \nonumber
\end{eqnarray}
The last inequality is because when $n$ is large, $x \geq n-k-1 > 2(k-2)$. Moreover, $\binom{x}{t}$ is an increasing function for $0<t<x/2$, so
when $x \leq n-4$, $\binom{x}{n-k-2}= \binom{x}{x-(n-k-2)} \leq \binom{x}{k-2}$.

Therefore we only need to verify \eqref{kk} for $n-k-1 \leq x< n-k$, which corresponds to $1 \leq m \leq n-k-1$.
For $m=1$, \eqref{kk} is obvious, so it suffices to look at the
case $m \geq 2$. The number of distinct $(k-1)$-subsets of $J_1$ or
$J_2$ is minimized when $|J_1 \cap J_2|=n-k-2$, which, by the 
inclusion-exclusion principle, gives $$N_2 \geq 2 \binom{n-k-1}{k-1} - \binom{n-k-2}{k-1} =
\binom{n-k-1}{k-1} + \binom{n-k-2}{k-2}.$$
So \eqref{counting} is also true for $2 \leq m \leq \binom{n-k-2}{k-2}+1$. It is easy to see that for $k \geq 3$ and $n$ sufficiently large,  $n-k-1 \leq \binom{n-k-2}{k-2}+1$.
For $k=2$, we have $x=n-3$ and \eqref{kk} becomes $\binom{n-3}{1} 
\geq \binom{n-3}{1} + \binom{n-3}{n-3} -1$, which is 
true and completes the proof. \qed

\noindent \textbf{Remark 1.} In order for all the inequalities to be correct, we only need $n>Ck^2$. By carefully analyzing the above computations, 
one can check that $C=500$ is enough.

\noindent \textbf{Remark 2.} Note that in the proof, the equality
\eqref{counting} holds in two different cases. The first case is when
$m=1$, which means $x_1 + x_{n-k+2}+ \cdots + x_n<0$ but any other
$k$-sums involving $x_1$ are nonnegative. All the other nonnegative
$k$-sums are formed by $x_2$ together with any $(k-1)$-subsets not
containing $x_{n-k+2}, \cdots, x_n$. This case is realizable by the
following construction: $x_1=k(k-1)n$, $x_2=n-2$, $x_3 = \cdots =
x_{n-k+1} = -1 $, $x_{n-k+2} = \cdots = x_n = -(kn+1)$. The second case
is in \eqref{kk} when $x=n-4$ and $x=n-k-1$ holds simultaneously, which
gives $k=3$. In this case, $m=\binom{n-3}{n-4}=n-3$, and the Kruskal-Katona
theorem holds with equality 
for the $(n-4)$-subsets $J_1, \cdots, J_{n-3}$. That is
to say, the negative $3$-sums using $x_1$ are $x_1 + x_i + x_n$ for
$3 \leq i \leq n-1$, while the nonnegative $3$-sums containing $x_2$
but not $x_1$ are $x_2+x_i + x_j$ for $3 \leq i<j \leq n-1$. This case
can also be achieved by setting $x_1=x_2=1$, $x_3= \cdots = x_{n-1}
= \frac{1}{2(n-3)}$, and $x_n=-\frac{3}{2}$. For large $n$, these are
the only two possible configurations achieving equality in Theorem
\ref{hm}.

Next we prove Theorem \ref{moderate}, which states that if 
$\sum_i x_i \geq 0$ and 
no $x_i$  is moderately large, then at least a constant proportion of the $\binom{n}{k}$ $k$-sums are nonnegative.\\
\\
\noindent \textbf{Proof of Theorem \ref{moderate}.}
Suppose $t$ is the largest
integer so that there are $t$ subsets $S_1, \cdots, S_t$ such that for any $1 \leq j \leq t$, $S_j$ is disjoint from $\{1, \cdots, j\}$, has at most $j(k-1)$ elements, and
$$x_1 + \cdots +x_j + \sum_{i \in S_j} x_i < 0.$$
By the maximality of $t$, the sum of $x_{t+1}$ and any $k-1$ numbers $x_i$ with indices from $[n] \backslash (\{1, \cdots, t+1\} \cup S_t)$ is nonnegative, so there are at least $\binom{n-tk-1}{k-1}$ nonnegative $k$-sums using $x_{t+1}$.
Since $x_{t+1}$ is not $(1-\delta)$-moderately large,
$$\binom{n-tk-1}{k-1} < (1-\delta) \binom{n-1}{k-1}.$$
For sufficiently large $n$, this is asymptotically equivalent to
$$\left (1-\dfrac{tk}{n} \right)^{k-1} < 1-\delta.$$
Since $$\left(1-\frac{tk}{n}\right)^{k-1} > 1-\frac{tk(k-1)}{n},$$
we have $$t > \dfrac{n}{k^2} \delta$$
Recall that by Corollary \ref{cor_numberofedges}, for each $i=1, \cdots, \frac{n}{k^2} \delta$, $x_i$ gives at least $\frac{1}{14} \binom{n-(i-1)k-1}{k-1}$ nonnegative $k$-sums, therefore
\begin{eqnarray}
N &\geq& \dfrac{1}{14} \left[ \binom{n-1}{k-1} + \cdots+ \binom{n-(\frac{n}{k^2}\delta)k-1}{k-1} \right] \nonumber \\
& \geq & \dfrac{n \delta }{14k^2} \left(1-\frac{\delta}{k}\right)^k \binom{n-1}{k-1} \nonumber \\
& > & \dfrac{\delta(1-\delta)}{14 k} \binom{n}{k}\nonumber
\end{eqnarray}
Setting $g(\delta, k) = \dfrac{\delta(1-\delta)}{14k}$ completes
the proof.
\qed

\section{Concluding remarks}
\label{concluding}
\begin{list}{\labelitemi}{\leftmargin=1em}

\item In this paper, we have proved that if $n>33k^2$,
any $n$ real numbers with a nonnegative sum have at least
$\binom{n-1}{k-1}$ nonnegative $k$-sums, thereby verifying the
Manickam-Mikl\'{o}s-Singhi conjecture in this range.  Because of the
inequality $\binom{n-2k}{k}+C\binom{n-1}{k-1}\geq \binom{n-1}{k-1}$ we
used, our method will not give a better range than the quadratic one,
and we did not try hard to compute the best constant in the quadratic
bound. It would be interesting to decide if 
the Manickam-Mikl\'{o}s-Singhi conjecture
can be verified for a linear range $n>ck$. Perhaps some algebraic methods
or structural analysis of the extremal configurations will help.

\item Feige \cite{feige} conjectures that the constant $1/13$ in Lemma
\ref{probabilistic_inequality} can be improved to $1/e$. This is a special case of a more general
question suggested by Samuels \cite{samuels_chebyshev}. He asked to determine, for a fixed $m$, the
infimum of $Pr(X_1+\cdots+X_k < m)$, where the infimum is taken over all possible collections of
$k$ independent nonnegative random variables $X_1, \cdots, X_k$ with given expectations
$\mu_1,\cdots,\mu_k$. For $k=1$ the answer is given by Markov's inequality. Samuels \cite{samuels_chebyshev, samuels2}
solved this question for $k \leq 4$, but for all $k \geq 5$ his problem is still completely
open.

\item Another intriguing objective is to prove the conjecture by Erd\H{o}s which states that
the maximum number of edges in an $r$-uniform hypergraph $H$ on $n$ vertices with matching number $\nu(H)$
is exactly \begin{equation*} \max \bigg \{\binom{r[\nu(H)+1]-1}{r},~
\binom{n}{r}-\binom{n-\nu(H)}{r} \bigg\}. \end{equation*}
The first number corresponds to a clique
and the second case is the complement of a clique. When $\nu(H)=1$, this conjecture is exactly the
Erd\H{o}s-Ko-Rado theorem \cite{erdos-ko-rado}. Erd\H{o}s also verified it for $n>c_r \nu(H)$ where
$c_r$ is a constant depending on $r$. Recall that in our graph $H$ we have $\nu(H) \leq n/k$ and
$r$ here is equal to $k-1$, so if Erd\H{o}s' conjecture is true in general, we can give a direct
proof of constant edge density in the complement of $H$. In this way we can avoid using fractional
matchings in our proof. But even without the application here, this conjecture is interesting in
its own right.  The fractional version of Erd\H{o}s' conjecture is also very interesting. In its
asymptotic form it says that if $H$ is an $r$-uniform $n$-vertex hypergraph with fractional
matching number $\nu^*(H)=xn$, where $0 \leq x <1/r$, then \begin{equation} \label{erdos} e(H) \leq
(1+o(1)) \max \big\{(rx)^r, 1-(1-x)^r \big \} {n \choose r}. \end{equation}

\item
As pointed out to us by Andrzej Ruci\'nski, part of our reasoning in
Section 3 implies that the function $A(n,k)$ defined in the first
page is precisely $\binom{n}{k}$ minus the maximum possible number
of edges in a $k$-uniform hypergraph on $n$ vertices with
fractional covering number strictly smaller than $n/k$. Indeed,
given $n$ reals $x_1, \ldots ,x_n$ with sum zero and only $A(n,k)$
nonnegative $k$-sums,  we may assume that the absolute value of
each $x_i$ is smaller than $1/k$ (otherwise simply multiply all of
them by a sufficiently
small positive real.) Next, add a sufficiently small  positive
$\epsilon$ to each $x_i$, keeping each $x_i$ smaller than $1/k$ and
keeping the sum of any
negative $k$-tuple below zero (this is clearly possible.) Note 
that the sum of these new reals, call them $x'_i$, is strictly
positive and the number  of positive $k$-sums is $A(n,k)$.
Put $\nu(i)=1/k-x'_i$ and observe that $\sum_i \nu (i) <n/k$ and
the $k$-uniform hypergraph  whose edges are all $k$-sets $e$ for
which $\sum_{i \in e} \nu(i) \geq 1$ has exactly
$\binom{n}{k}-A(n,k)$ edges.   
Therefore, there is a $k$-uniform  hypergraph
on $n$ vertices with fractional covering number strictly smaller
than $n/k$ and at least $\binom{n}{k}-A(n,k)$ edges.
Conversely, given a $k$-uniform hypergraph $H$ on $n$ vertices 
and a fractional covering of it
$\nu:V(H) \mapsto [0,1]$ with $\sum_i \nu (i)=n/k-\delta <n/k$ 
and $\sum_{i
\in e} \nu (i) \geq 1$ for each $e \in E(H)$, one can define
$x_i=\frac{1}{k}-\frac{\delta}{n}-\nu(i)$  to get a set of $n$ 
reals whose  sum  is zero, in which the number of nonnegative
$k$-sums is at most $\binom{n}{k} -|E(H)|$ (as the sum of the
numbers $x_i$ for every
$k$-set forming an edge of $H$ is  at most $1-\frac{k \delta}{n}
-1<0$). This implies the desired equality, showing that the problem
of determining $A(n,k)$ is equivalent to that of finding the
maximum possible number of edges of a $k$-uniform hypergraph on $n$
vertices with fractional covering number strictly smaller than
$n/k$. Note that this is equivalent to the problem of settling the
fractional version of the conjecture of Erd\H{o}s for the extremal
case of fractional matching number $<n/k$.

\item Although the fractional version of Erd\H{o}s' conjecture is still widely open in general,
our techniques can be used to make some progress on this problem. Combining the approach from Section 3
and the above mentioned results of Samuels we verified, in joint work with Frankl, R\"odl and Ruci\'nski
\cite{AFHRRS}, conjecture (\ref{erdos}) for certain ranges of $x$ for $3$ and $4$-uniform
hypergraphs. These results can be used to study a Dirac-type question of Daykin and H\"aggkvist
\cite{DH} and of H\'an, Person and Schacht \cite{HPS} about perfect matchings in hypergraphs.

For an $r$-uniform hypergraph $H$ and for $1 \leq d \leq r$, let $\delta_d(H)$ denote the minimum
number of edges containing a subset of $d$ vertices of $H$, where the minimum is taken over all
such subsets. In particular, $\delta_1(H)$ is the minimum vertex degree of $H$. For $r$ that
divides $n$, let $m_d(r,n)$ denote the smallest number, so that any $r$-uniform hypergraph $H$ on
$n$ vertices with $\delta_d(H) \geq m_d(r,n)$ contains a perfect matching. Similarly, let
$m^*_d(r,n)$ denote the smallest number, so that any $r$-uniform hypergraph $H$ on $n$ vertices
with $\delta_d(H) \geq m^*_d(r,n)$ contains a perfect fractional matching.

Together with Frankl, R\"odl and Ruci\'nski \cite{AFHRRS} we proved that for all $d$ and $r$,
$m_d(r,n) \sim m^*_d(r,n)$ and further reduced the problem of determining the asymptotic behavior
of these numbers to some special cases of conjecture (\ref{erdos}). Using this relation we were
able to determine $m_d(r,n)$ asymptotically for several values of $d$ and $r$, which have not been known
before. Moreover, our approach may lead to a solution of the general case as well, see \cite{AFHRRS} for the details.

\end{list}
\vspace{0.2cm}

\noindent
{\bf Acknowledgment} We would like to thank Andrzej Ruci\'nski for 
helpful discussions and comments, and Nati Linial for inspiring
conversations and intriguing questions which led us to the 
results in Section 4.

\end{document}